\documentclass[12pt]{amsart}
\usepackage{amssymb}

\newtheorem{Thm}{Theorem}[section]

\newtheorem{Lem}[Thm]{Lemma}
\newtheorem{Prop}[Thm]{Proposition}

\theoremstyle{definition}

\theoremstyle{remark}

\def \LIM{\text{LIM }}

\begin{document}

\title{Solution of a problem of Peller concerning similarity
}
\author{N. J. Kalton}
\address{Department of Mathematics \\
University of Missouri-Columbia \\
Columbia, MO 65211
}
\thanks{The first author was supported by NSF grant DMS-9870027}

\email[N. J. Kalton]{nigel@math.missouri.edu}

\author {C. Le Merdy}
\address{Equipe de Math\'ematiques - UMR 6623, Universit\'e de
Franche-Comt\'e, F-25030 Besan\c con cedex}
\email[C. Le Merdy]{lemerdy@math.univ-fcomte.fr}


\begin{abstract}We answer a question of Peller by showing that for any
$c>1$ there exists a power-bounded  operator $T$ on a Hilbert
space with the property that any
operator
$S$ similar to
$T$ satisfies
$\sup_n\|S^n\|>c.$
\end{abstract}

\maketitle

\section{Introduction} In this note we answer a question due
to Peller \cite{Pe} which has also recently been raised by Pisier
\cite{Pi} p.114.
 Peller's question is whether, for any
$\epsilon>0$, every power-bounded operator $T$ is similar to an operator
$S$ with $\sup_n\|S^n\|< 1+\epsilon.$

 It was shown by Foguel \cite {Fo} in 1964 that there
is a power-bounded operator $T$ on a Hilbert space ${{\mathcal H}}$  which is not
is not similar to a contraction.  It was later shown by Lebow that this
example
is not polynomially bounded \cite{L}; for other examples see \cite{Bo}
and
\cite{Pi}, Chapter 2.  Recently Pisier \cite{Pi2} answered a problem
raised by Halmos  by constructing an
operator which is polynomially bounded and not similar to a contraction.

We shall construct a family of counter-examples to Peller's question.
These counter-examples have a rather simple structure.  Let $w$ be an
$A_2-$weight on the circle $\mathbb T$ and let ${H}^2(w)$ be the closed
linear span of $\{e^{in\theta}:\ n\ge 0\}$ in $L^2(w).$  We consider an
operator
$$
T(\sum_{n=0}^{\infty}a_ne^{in\theta})=\sum_{n=0}^{\infty}\lambda_na_ne^{in\theta}
$$
where $(\lambda_n)_{n=0}^{\infty}$ is a monotone increasing sequence of
positive reals with $\lambda_n\uparrow 1$ and $\lambda_n<1$ with
$$ \lim_{n\to\infty}\frac{1-\lambda_{n+1}}{1-\lambda_n}=0.$$
 For such operators we can prove a rather precise result (Theorem
\ref{main}):
\begin{equation}\label{eq}
 \inf\{\sup_n\|(A^{-1}TA)^n\|: \ A \text{
invertible}\}=\sec(\frac{\pi}{2p})\end{equation}
where $p=\sup\{a:w^a\in A_2\}.$  By taking simple choices of
$A_2-$weights where $p<\infty$ we can create a family of
counter-examples.

The proof of Theorem \ref{main} depends heavily on estimates for the norm
of the Riesz projection in Section 2 particularly Theorem \ref{opt}.
These results can be obtained by a careful reading of the classical work
of Helson and Szeg\"o \cite{HS}  on $A_2-$weights (cf. \cite{G}).
However,
we present a self-contained argument, in which the reader will recognize
many similarities with the Helson-Szeg\"o theory.

We also show that our examples can only be polynomially bounded in the
trivial situation when $w$ is equivalent to the constant function and
then $T$ is similar to contraction. We also note that the case $p=\infty$
in (\ref{eq}) (when Peller's conjecture holds for $T$) corresponds to the
case when $\log w$ is in the closure of
$L^{\infty}(\mathbb T)$
in $BMO(\mathbb T).$

\section{The norm of the Riesz projection on weighted $L^2-$spaces}
\setcounter{equation}{0}

We start by recalling an easy lemma concerning projections on a Hilbert
space.

\begin{Lem}\label{elem} Let $E$ and $F$ be closed subspaces of a Hilbert
space
${{\mathcal H}}$
so that $E+F$ is dense in ${{\mathcal H}}.$  Suppose $0\le \varphi<\pi/2.$  In order
that there is a projection $P$ of ${{\mathcal H}}$ onto $E$ with $F=\text{ker }P$
with $\|P\|\le \sec\varphi$ it is necessary and sufficient that
$$ |(e,f)|\le \sin\varphi \|e\|\|f\| \qquad e\in E,\ f\in F.$$
\end{Lem}

{\it Remark.}  Note that a consequence of Lemma \ref{elem} is that if $P$
is any non-trivial projection on a Hilbert space then $\|P\|=\|I-P\|.$

Now let $\mathbb T$ be the unit circle (which we identify with
$(-\pi,\pi]$ in the usual way) equipped with the standard Haar measure
$d\theta/2\pi.$  Let $\mu$ be any finite positive Borel measure on
$\mathbb T.$  We denote by $L^2(\mu)=L^2(\mathbb T;\mu)$ the
corresponding
weighted $L^2-$space; if $\mu$ is absolutely continuous with respect to
Haar measure so that $d\mu=(2\pi)^{-1}w(\theta)d\theta$ then we write
$L^2(w).$  We refer to any nonnegative $w\in L^1(\mathbb T)$
so that $w>0$ on a set of positive measure as a weight.

Suppose $w$ is a weight. We recall that
${H}^2(w)$ is the closed subspace of
$L^2(w)$ generated by the functions $\{e^{in\theta}:\ n\ge 0\}.$
We recall that $w$ is an
 {\it $A_2-$weight} if there is a bounded
projection $R$ of $L^2(w)$ onto ${H}^2(\mu)$ with $R(e^{in\theta})=0$
if
$n<0.$  In this case we always have that $w>0$ a.e., $w^{-1}$ is an
$A_2-$weight and
$L^2(w)\subset L^1;$ the operator $R$ must coincide with the Riesz
projection $Rf\sim\sum_{n\ge 0}\hat f(n)e^{in\theta}.$  Let us denote by
$\|R\|_w$ the norm of the Riesz projection on $L^2(w).$  Note that for an
$A_2-$weight ${H}^2(w)={H}^1\cap L^2(w)$.  In particular we can define
$f(z)=\sum_{n\ge 0}\hat f(n)z^n$ for $|z|<1$.

The following Proposition can be derived from the classical
work of Helson-Szeg\"o \cite{HS} or \cite{G}.  However, we give a
self-contained direct proof.

\begin{Prop}\label{helson} Let $w$ be a weight function on $\mathbb T.$
Assume $0\le \varphi<\frac{\pi}{2}.$ The following conditions are
equivalent:
\newline
(1) $w$ is an $A_2-$weight and $\|R\|_w \le \sec\varphi.$\newline
(2) There exists $h\in H^1$ so that $|w-h|\le w\sin\varphi $ a.e.
\end{Prop}

\begin{proof}  First note that by Lemma \ref{elem} (1) is equivalent to
\begin{equation} \label{elem2}
\left|\int_{-\pi}^{\pi}f(\theta)g(\theta)w(\theta)\frac{d\theta}{2\pi}
\right| \le \sin\varphi
\left(\int_{-\pi}^{\pi}|f(\theta)|^2w(\theta)\frac{d\theta}{2\pi}\right)^{1/2}
\left(\int_{-\pi}^{\pi}|g(\theta)|^2w(\theta)\frac{d\theta}{2\pi}\right)^{1/2}
\end{equation},
whenever $f,g\in H^2(w)$ with $g(0)=0.$

To prove (1) implies (2) we note that if $w$ is an $A_2-$weight so that
$\log w\in L^1$ we can find an outer function $F\in H^2$ so that
$w=|F|^2$ a.e..  Then (\ref{elem2}) gives
$$\left| \int_{-\pi}^{\pi}fgwF^{-2}\frac
{d\theta}{2\pi}\right| \le \sin\varphi
\left(\int_{-\pi}^{\pi}|f|^2\frac{d\theta}{2\pi}\right)^{1/2}
\left(\int_{-\pi}^{\pi}|g|^2\frac{d\theta}{2\pi}\right)^{1/2},$$
for $f,g\in H^2$ with $g(0)=0.$
 This in turn implies that
$$ \left|\int_{-\pi}^{\pi}fwF^{-2}\frac{d\theta}{2\pi}\right| \le
\sin\varphi \|f\|_1$$
for all $f\in H^1,$ with $f(0)=0.$ By the Hahn-Banach Theorem this
implies there
exists $G\in H^{\infty}$ so that $\|wF^{-2}-G\|_{\infty}\le \sin\varphi$
or $|w-h|\le w\sin\varphi$ where $h=F^2G\in H^1.$

For the reverse direction just note that if $f,g\in H^2(w)$ with
$g(0)=0$ then
$$ \int_{-\pi}^{\pi} fgw \frac{d\theta}{2\pi} =\int_{-\pi}^{\pi}
fg(w-h)\frac{d\theta}{2\pi}$$
so that (\ref{elem2}) follows from the Cauchy-Schwartz inequality.
\end{proof}

Let us isolate a simple special case of the above proposition.

\begin{Prop}\label{new}  Let $0\neq f \in H^1$ be such that $\arg
f(\theta)\le \varphi<\frac{\pi}{2} $almost everywhere.   If $f$ is
not
identally zero then $w=\Re f$ is an $A_2-$weight for which $\|R\|_w \le
\sec\varphi.$\end{Prop}

\begin{proof} In this case $w=\Re f\ge 0$ a.e. and $|\Im f| \le
\tan\varphi w$ a.e. Furthermore:
$$ |w-\cos^2\varphi f|^2 \le (\sin^4\varphi + \cos^4\varphi\tan^2
\varphi)w^2 \le \sin^2\varphi w^2$$
a.e., so that we obain the result from Proposition
\ref{helson}.\end{proof}

{\it Remark.}  Suppose $0<\alpha<1$ and $f\in H^1(\mathbb D)$ is given by
$$ f(z)=\frac{(z-1)^{\alpha}}{(z+1)^{\alpha}}$$
(taking the usual branch of $w\mapsto w^{\alpha}.$)  Then $$w=\Re
f=\cos\frac{\alpha\pi}{2}\tan^{\alpha}\frac{\theta}{2}.$$ It follows
that
\begin{equation}\label{weight}
 \|R\|_{|\tan^{\alpha}(\theta/2)|}\le \sec\frac{\alpha\pi}{2}.
 \end{equation}
In fact (\ref{weight}) is well-known (see \cite{Krup}, for
example).  We are grateful to Igor Verbitsky for bringing this reference
to our        attention.

We will say that two weights $v,w$ are equivalent ($v\sim w$) if
$v/w,w/v\in L^{\infty}.$

\begin{Thm}\label{opt} Suppose $w$ is an $A_2-$weight on $\mathbb T.$
Then
$$ \inf\{\|R\|_v:\ v\sim w\}= \sec(\frac{\pi}{2p})$$
where
$$ p= \sup\{ a>0:  \ w^a\in A_2\}.$$
\end{Thm}

\begin{proof} First suppose $v\sim w$ and $\|R\|_v=\sec\psi$ where $0\le
\psi<\pi/2.$  Then there exists $h\in H^1$ with $|v-h|\le
v\sin\psi$ a.e.  In particular, $|\arg h|\le \psi$ a.e. and so $h$ maps
$\mathbb D$ into the same sector.  It follows that we can define $h^r\in
H^{1/r}$ for all $r>0.$  Choose $r$ so that $r\psi<\pi/2,$ and let
$g=h^r.$ Then
$\Re
g\ge 0$ and $|\Im g| \le \tan (r\psi) \Re g$ so that $g\in
H^1.$  Now by Proposition \ref{new} we have that
 $\Re g$ is an
$A_2-$weight.
However $\Re g\sim |h|^r\sim w^r$ so that $r\le p.$  We deduce that
$\psi\ge \pi/(2p).$

For the converse direction assume that $w^r$ is an $A_2-$weight. Then
there exists
$h\in H^1$ so that $|w^r-h|\le w^r\sin\psi$ where $0\le \psi<\pi/2.$
Arguing as above we have $g=h^{1/r}\in H^1$ and $\Re g$ is an
$A_2-$weight with $\|R\|_{\Re g}\le \sec(\psi/r).$  Note that $\Re g\sim
w$, and this establishes the other direction.\end{proof}

{\it Remark.}  If we now let
$w(\theta)=|\tan\frac{\theta}{2}|^{\alpha}$ where
$0<\alpha<1$ then we can apply (\ref{weight}) to deduce that, for this
particular weight the infimum is attained, i.e.
\begin{equation}\label{weight2} \inf\{\|R\|_v:\ v\sim w
\}=\|R\|_{|\tan^{\alpha}(\theta/2)|}=\sec\frac{\alpha
\pi}{2}.\end{equation}

\section{Multipliers} \label{mult}  \setcounter{equation}{0}

Suppose $(e_n)_{n=0}^{\infty}$ be any Schauder basis of a Hilbert
space
${\mathcal H}$; note that
we do not assume $(e_n)$ to be orthonormal or even unconditional. Let
$(P_n)$ be the associated partial sum operators
$P_n(\sum_{k=0}^{\infty}a_ke_k)=\sum_{k=0}^na_ke_k.$  Let $Q_n=I-P_n$ and
note
that $\|Q_n\|=\|P_n\|$ for all $n\ge 0.$  Since $(e_n)$ is a basis we
have that $\sup_n\|P_n\|=b<\infty$ where $b$ is the {\it basis
constant}.
 We
call an operator $T:{\mathcal H}\to {\mathcal H}$ a {\it monotone multiplier}
(with respect to
the given basis) if there is an increasing sequence
$(\lambda_k)_{k=0}^{\infty}$ in $\mathbb R$ so that
$0\le\lambda_k\le  1$ so that
$$ T(\sum_{k=0}^{\infty}a_ke_k)=\sum_{k=0}^{\infty}\lambda_ka_ke_k.$$

\begin{Lem} \label{bv}  If $T$ is defined as above then $T$ is
(well-defined and) bounded and $\sup_n\|T^n\|\le b.$\end{Lem}

\begin{proof} It is enough to show $T$ is bounded and $\|T\|\le b$ since
$T^n$ is also a monotone multiplier.  To see this note that if
$(a_k)_{k=0}^{\infty}$ is finitely nonzero  and
$x=\sum_{k=0}^{\infty}a_ke_k,$ then
$$ Tx = \lambda_0x+\sum_{k=1}^{\infty}(\lambda_k-\lambda_{k-1})Q_kx$$
so that $\|Tx\| \le \sup_n\|Q_n\|=b.$\end{proof}

We shall say that $T$ is a {\it fast monotone multiplier} if in addition,
$\lambda_k<1$ for all $k$ and
\begin{equation}\label{fmm}
\lim_{k\to\infty}\frac{1-\lambda_k}{1-\lambda_{k-1}}=0.\end{equation}

\begin{Lem} \label{fmm2}  Suppose $T$ is a fast monotone multiplier.
Then there is an increasing sequence of integers $(N_n)_{n=0}^{\infty}$
so that $\lim_{n\to\infty}\|T^{N_n}-Q_n\|=0.$\end{Lem}

\begin{proof} Note that if $x=\sum_{k=0}^{\infty}a_ke_k$ then
$$ T^{N_n}x-Q_nx = \sum_{k=0}^{n}\lambda_k^{N_n}a_ke_k
-(1-\lambda_{n+1}^{N_n})Q_nx
+\sum_{k=n+1}^{\infty}(\lambda_k^{N_n}-\lambda_{n+1}^{N_n})a_ke_k$$
whence a calculation as in Lemma \ref{bv} gives
$$ \|T^{N_n}x-Q_nx\| \le b\lambda_n^{N_n}\|P_nx\|
+(b+1)(1-\lambda_{n+1}^{N_n})
\|Q_nx\|.$$
It follows that
$$ \|T^{N_n}-Q_n\| \le
b(b\lambda_n^{N_n}+(b+1)(1-\lambda_{n+1}^{N_n})).$$
It remains therefore only to select $N_n$ so that
$\lim_{n\to\infty}\lambda_n^{N_n}=0$ and
$\lim_{n\to\infty}\lambda_{n+1}^{N_n}=1.$

For convenience we write $\lambda_n=e^{-\nu_n}$ where
$\nu_{n}/\nu_{n+1}=\kappa_n^2$ and $\kappa_n\to\infty.$  For any
$n\ge 0,$ pick
$N_n$
to be the greatest integer so that $N_n\nu_n^{1/2}\nu_{n+1}^{1/2}\le
1.$ Then
$$ N_n\nu_{n-1}^{1/2}\nu_{n}^{1/2}\ge \frac{N_n}{N_n+1}$$
and $\lim N_n=\infty.$

 Now
$$ N_{n}\nu_{n} \ge\frac{N_{n} \kappa_n}{N_n+1}$$
and
$$ N_n\nu_{n+1} \le \kappa_n^{-1}.$$
This yields the desired result.\end {proof}

We now turn to the case when ${\mathcal H}=H^2(w)$ where $w$ is an $A_2-$weight and
$e_k(\theta)=e^{ik\theta}$ for $k\ge 0.$

\begin{Lem}\label{basis}  The basis constant of $(e_k)_{k=0}^{\infty}$ in
$H^2(w)$ is given by $b=\|R\|_w.$\end{Lem}

\begin{proof} In fact $Q_{n-1}f=e_nR(e_{-n}f)$ so it
is clear that $\|Q_{n-1}\|\le \|R\|_w.$   For the other direction suppose
$f$ is a
trigonometric polynomial in $L^2(w).$  Then for large enough $n$ we have
$e_nf\in H^2(w)$ and then $Rf= e_{-n}Q_{n-1}(e_nf)$.  This quickly yields
$\|R\|_w\le b.$\end{proof}

\begin{Thm}\label{main} Let $w$ be an $A_2-$weight on $\mathbb T$ and let
$T:H^2(w)\to H^2(w)$ be a fast monotone multiplier corresponding
to the sequence $(\lambda_n)$. Then
\begin{equation}\label{formula} \inf\{\sup_n\|(A^{-1}TA)^n\|: \ A \text{
invertible}\}
=
\sec{\frac{\pi}{2p}}\end{equation}
where
$$ p= \sup\{ a>0:  \ w^a\in A_2\}.$$
\end{Thm}

\begin{proof}
 We shall prove that if $\sigma\ge 1$ then the existence of
an invertible $A$ so that $\sup_n\|(A^{-1}TA)^n\|\le \sigma$ is
equivalent to the existence of a weight $v$ equivalent to $w$ so that
$\|R\|_v\le \sigma.$  Once this is done, the result follows from  Theorem
\ref{opt}.

In one direction this is easy.  Assume $v$ equivalent to $w$ and
$\|R\|_v\le \sigma.$  This means that there is an equivalent
inner-product norm on $H^2(w)$ in which the basis constant of
$(e_k)_{k=0}^{\infty}$ bounded by $\sigma.$  It follows from Lemma
\ref{bv} that in this equivalent norm we have $\sup_n\|T^n\|_v\le
\sigma.$  Hence $T$ is similar to an operator $A^{-1}TA$ such that
$\sup\|(A^{-1}TA)^n\|\le \sigma.$

We now consider the converse.  Let $S:H^2(w)\to H^2(w)$ be the operator
$Sf=e_1f.$  Suppose $A$ is an invertible operator such that
$\|(A^{-1}TA)^n\|\le \sigma.$  We will define a new inner-product on
$H^2(w)$ by
$$ \langle f,g\rangle = \LIM (AS^nf,AS^ng)$$
where $\LIM$ denotes any Banach limit (see e.g. \cite{Co} p. 85). Since $S$ is an isometry on
$H^2(w)$ and $A$ is invertible this defines an equivalent inner-product
$|\cdot|$
norm on $H^2(w).$  Now for any $f\in H^2(w)$ and fixed $m\in\mathbb N$
we have
$$ \lim_{n\to\infty}\|AQ_{m+n}S^nf-AT^{N_{m+n}}S^nf\|=0$$
where $(N_n)$ is given in Lemma \ref{fmm2}.  Hence
$$\limsup_{n\to\infty}(\|AQ_{m+n}S^nf\|^2-\sigma^2 \|AS^nf\|^2)\le 0.$$
Now
$$ |Q_mf|^2=\LIM \|AS^nQ_mf\|^2 =\LIM\|AQ_{m+n}S^nf\|^2 \le \sigma^2
|f|^2.$$
Thus with respect to the new norm $|\cdot|$ the basis constant is at most
$\sigma.$

Now let $c_k=\langle e_0,e_k\rangle$ for $k\ge 0$ and let
$c_k=\overline{c}_{-k}$ when $k<0.$ Then
it follows easily that $\langle e_k,e_l\rangle = c_{l-k}$ for all $k,l$
and that for all finitely nonzero sequences $(a_k)$ of complex numbers we
have that
$$ \sum_{k,l}a_k\overline{a}_lc_{k-l}\ge 0.$$
This implies (see \cite{K} p. 38) that there is a finite positive measure
$\mu$ on
$\mathbb T$ so that
$$ \int e^{-ik\theta}d\mu(\theta)=c_k.$$
Thus
$$ \langle f,g\rangle =\int f\overline{g}\,d\mu.$$
However this norm is equivalent to the original norm so that $\mu$ is
absolutely continuous with respect to Lebesgue measure and of the form
$(2\pi)^{-1}v(\theta)d\theta$ where $v\sim w.$

It follows that in $H^2(v)$ the basis constant of the exponential basis
is at most $\sigma$ and so by Lemma \ref{basis} we have $\|R\|_v\le
\sigma$ and the proof is complete.\end{proof}

We can now give explicit examples by taking the weights
$w(t)=|\theta|^{\alpha}$ where $0<\alpha<1.$  It is clear that in
Theorem
\ref{main} we have $p=\alpha^{-1}$ and so for any fast monotone
multiplier we have
$$ \inf\{\sup_n\|(A^{-1}TA)^n\|: \ A \text{ invertible}\} =
\sec{\frac{\pi\alpha}{2}}>1.$$  Note that we are essentially using here
the original example of a conditional basis for Hilbert space due to
Babenko \cite{B}.  We can also utilize (\ref{weight2}) to show that for
this example the infimum in (\ref{formula}) is actually attained.  In
general the infimum in (\ref{formula}) need not be attained; this is
will be seen easily from Theorem \ref{last} below.

\begin{Thm}\label{contraction}  Let $w$ be an $A_2-$weight and suppose
$T:H^2(w)\to H^2(w)$ is a fast monotone multiplier, corresponding to the
sequence $(\lambda_n).$
 Then the following are equivalent: \newline
(i) $T$ is similar to a contraction.\newline
(ii) $T$ is polynomially bounded. \newline
(iii) $w\sim 1.$
\end{Thm}

\begin{proof} That (i) implies (ii) is a consequence of von Neumann's
inequality (see \cite{Pi}). Similarly (iii) implies (i) is trivial.  It
therefore remains to prove that (ii) implies (iii). We shall treat the
case when the $\lambda_k$ are distinct; small modifications are necessary
in the other cases. We shall also suppose the measure
$d\mu=(2\pi)^{-1}w(\theta)d\theta$ is a probability measure so that
$\|e_k\|=1$ for all $k.$

 First note that if $f\in H^{\infty}(\mathbb D)$ then for any $r<1,$
then $f_r(T)$ is well-defined where $f_r(z)=f(rz)$ and if  $T$ is
polynomially bounded we have an estimate $$\|f_r(T)\|\le
C\|f\|_{H^{\infty}(\mathbb D)},$$ or equivalently
$$ \|\sum_{k=0}^{\infty}f(r\lambda_k)a_k e_k\| \le
C\|f\|_{H^{\infty}(\mathbb D)} \|\sum_{k=0}^{\infty}a_ke_k\|$$
whenever $(a_k)$ is finitely non-zero. Letting $r\to 1$ we obtain
$$ \|\sum_{k=0}^{\infty}f(\lambda_k)a_k e_k\| \le
C\|f\|_{H^{\infty}(\mathbb D)} \|\sum_{k=0}^{\infty}a_ke_k\|$$

 Recall that by Carleson's theorem \cite{C} the sequence
$(\lambda_n)$ is
{\it interpolating}
(cf. \cite{G} p. 287-288) so that there is a constant $B$ such that for
any sequence
$\epsilon_k=\pm1$ there exists $f\in H^{\infty}(\mathbb D)$ with
$\|f\|_{H^{\infty}(\mathbb D)}\le B$ and $f(\lambda_k)=\epsilon_k$ for
all $k\ge 0.$   Hence
$$ \|\sum_{k=0}^{\infty}\epsilon_ka_k e_k\| \le
BC \|\sum_{k=0}^{\infty}a_ke_k\|$$
for all finitely non-zero sequences $(a_k).$   Hence
by the parallelogram law we have
$$ (BC)^{-1}(\sum_{k=0}^{\infty}|a_k|^2)^{1/2} \le
\|\sum_{k=0}^{\infty}a_ke_k\| \le BC (\sum_{k=0}^{\infty}|a_k|^2)^{1/2}$$
from which it follows that $w\sim 1.$  \end{proof}

We conclude by considering the cases when
$$ \inf\{\sup_n\|(A^{-1}TA)^n\|: \ A \text{ invertible}\} =  1.$$

\begin{Thm}\label{last}
  Let $w$ be an $A_2-$weight and suppose
$T:H^2(w)\to H^2(w)$ is a fast monotone multiplier, corresponding to the
sequence $(\lambda_n).$
 Then the following are equivalent: \newline
 (i) For any $\epsilon>0,$ $T$ is similar to an operator $S$ with $\sup_n
\|S^n\|<1+\epsilon.$
  \newline
 (ii) $\log w$ is in the closure of $L^{\infty}$ in $BMO.$  \newline
 (iii) $w^a\in A_2$ for every $a>0.$
 \end{Thm}

 \begin{proof}  The equivalence of (i) and (iii) is proved in Theorem
\ref{main}. The equivalence of (ii) and (iii) is due to Garnett and Jones
\cite{GJ} or \cite{G}, Corollary 6.6 and its proof (p.258-9).\end{proof}

\end{document}